\documentclass{article}
\usepackage{graphicx}         
\usepackage{amsmath}          
\usepackage{amsfonts}
\usepackage{fancyhdr}
\usepackage{a4}
\usepackage{float}
\def\a{\alpha}

\def\b{\beta}

\def\bC{\mathbf{C}}

\def\d{\delta}
\def\bd{\boldsymbol{\delta}}

\def\ds{\displaystyle}

\def\E{\mathbf{E}}
\def\e{\mathbf{e}}

\def\f{\mathbf{f}}

\def\G{\mathbf{G}}
\def\bg{\boldsymbol{\gamma}}

\def\I{\mathbf{I}}

\def\l{\lambda}

\def\norm#1{\parallel#1\parallel}

\def\O{\Omega}

\def\P{\mathbf{P}}
\def\phi{\varphi}

\def\bPhi{\mathbf{\Phi}}

\def\R{\mathbf{R}}
\def\r{\mathbf{r}}

\def\s{\mathbf{s}}

\def\bt{\boldsymbol{\tau}}

\begin{document}
\title{A pseudo active kinematic constraint for a biological living soft tissue: an effect of the collagen network}
\author{
Christian Bourdarias\thanks{Laboratoire de Math\'ematiques, Universit\'e de Savoie,
73376 Le Bourget du Lac, France, 
email:Christian Bourdarias @univ-savoie.fr}, 
St\'{e}phane Gerbi\thanks{Laboratoire de Math\'ematiques, Universit\'e de Savoie,
73376 Le Bourget du Lac, France, 
email:Stephane.Gerbi@univ-savoie.fr} ~and
Jacques Ohayon\thanks{ Universit\'e Joseph Fourier, Laboratoire TIMC-IMAG, Grenoble, France and 
Institut de l'Ing\'enierie et de l'Information de Sant\'e (IN3S)
Facult\'e de M\'edecine de Grenoble
Bâtiment Taillefer (Bureau B57)
38706 La Tronche Cedex, France,
email: Jacques.Ohayon@imag.fr}
}
\date{}
\maketitle
%
\begin{abstract}
Recent studies in mammalian hearts show that left ventricular wall thickening is an important mechanism for 
systolic ejection and that during contraction the cardiac muscle develops
significant stresses in the muscular cross-fiber direction. We suggested that the collagen
network surrounding the muscular fibers could account for these mechanical behaviors. 
To test this hypothesis we develop a  model for large deformation response of active, incompressible, nonlinear elastic and
transversely isotropic living soft tissue (such as cardiac or arteries tissues) in which we include a coupling effect
between the connective tissue and the muscular fibers. 
Then, a three-dimensional finite element formulation including this internal pseudo-active kinematic
constraint is derived. Analytical and finite element solutions are in a very good agreement. The numerical results show this wall thickening effect with an order of magnitude compatible with the experimental observations.
\end{abstract}

\section{Introduction}

It is known that the transverse shear along myocardial cleavage planes provides a mechanism for
a normal systolic wall thickening \cite{LTC95}. Indirect evidences indicate that the characteristics
of the passive extracellular connective tissue in the myocardium is an important determinant
of ventricular function (\cite{LY98}, \cite{R01}, \cite{CHWWGM962}). An appropriate constitutive law for the
myocardium should therefore incorporate the most important features of its microstructure. A sound
theoretical formulation for material laws of the active myocardium is essential for an accurate
mechanical analysis of the stresses in the ventricular wall during the whole cardiac cycle.
The wall stress distribution is one of the main factors governing the myocardial energetic
\cite{SBWCSM58}, the coronary blood flow \cite{CTMOL90}, the
cardiac hypertrophy \cite{R01}, and the fetal heart growth \cite{OCJUA01}.
To date we do not have any reliable technique to evaluate the stress in the cardiac
muscle, therefore, mechanical models are useful in cardiology to assess the functional capacities of
the human heart. Several numerical models using a finite element (FE) analysis have been
performed to simulate the left ventricular performance  (\cite{{HCAH91}}, \cite{UHM01}). The mechanical
behavior of the connective tissue is often assumed isotropic \cite{OC88}. This last assumption is not in
agreement with the experimental results obtained on a sample of active myocardial rabbit tissue. Lin
and Yin \cite{LY98} showed that, during an active equibiaxial stretch test, there are significant
stresses developed in the cross-fiber direction (more than 40\% of those in the fiber direction)
that cannot be attributed to nonparallel muscle fibers (MF).

Therefore, the purposes of this paper are to: (i) suggest a realistic pseudo-active kinematic
law coupling the passive connective tissue to the muscle fibers, which may explain the
developed tension in the cross-fiber direction observed by Lin and Yin (\cite{LY98}), (ii) formulate
an active three-dimensional material law for a nonlinear hyperelastic and incompressible
continuum medium, which takes care of these coupling effects, (iii) derive the related
three-dimensional finite element (FE) formulation, and (iv) test the accuracy of the
proposed numerical method.


\section{Microstructure of the cardiac tissue}

\subsection{Muscle fiber organization}

Anatomical observations have shown that the cardiac muscle tissue has a highly specialized architecture
\cite{S79}. This structure is composed primarily of cardiac muscle cells, or myocytes, that are 80 to
100 $\mu$m in length and are roughly cylindrical with cross-sectional dimensions of 10 to 20 $\mu$m.
These cells are arranged in a more or less parallel weave that we idealize as ``muscle fibers" (MF). We
shall denote the local direction of this group of cells by the unit vector $\f$ and refer to it also
as the local ``fiber" direction with the understanding that individual continuous MF do
not really exist. Experimental measurements have shown that the MF direction field defines
paths on a nested family or toroidal surfaces of revolution in the wall of the heart \cite{S79}.
These results show a continuously changing orientation $\f$ of the MF through the wall,
circumferential near the midwall and progressively more inclined with respect to the equatorial
plane when moving toward either the epicardium or the endocardium.

\subsection{The cardiac connective tissue organization}

Myocytes and coronary blood vessels are embedded in a complex extracellular matrix which
consists of collagen and elastin, mainly. Caulfield and Janicki  \cite{CJ97} used the scanning
electron microscope (SEM) to reveal the basic organization of this connective tissue network.
Their studies on the connective tissue of mammalian heart muscle give the description of the
extracellular structures and their arrangement relative to cardiac muscle cells.
They described the three following classes of connective tissue
organization: (i) interconnections between myocytes, (ii) connections between myocytes and
capillaries and, (iii) a collagen weave surrounding group of myocytes. When viewed by SEM, groups of
myocytes can be seen to be encompassed by a rather prominent meshwork of fibrillard collagen, and
short collagen struts attach the myocytes subjacent to this meshwork to it.

\begin{figure}[H]
\centering\includegraphics[scale=0.2]{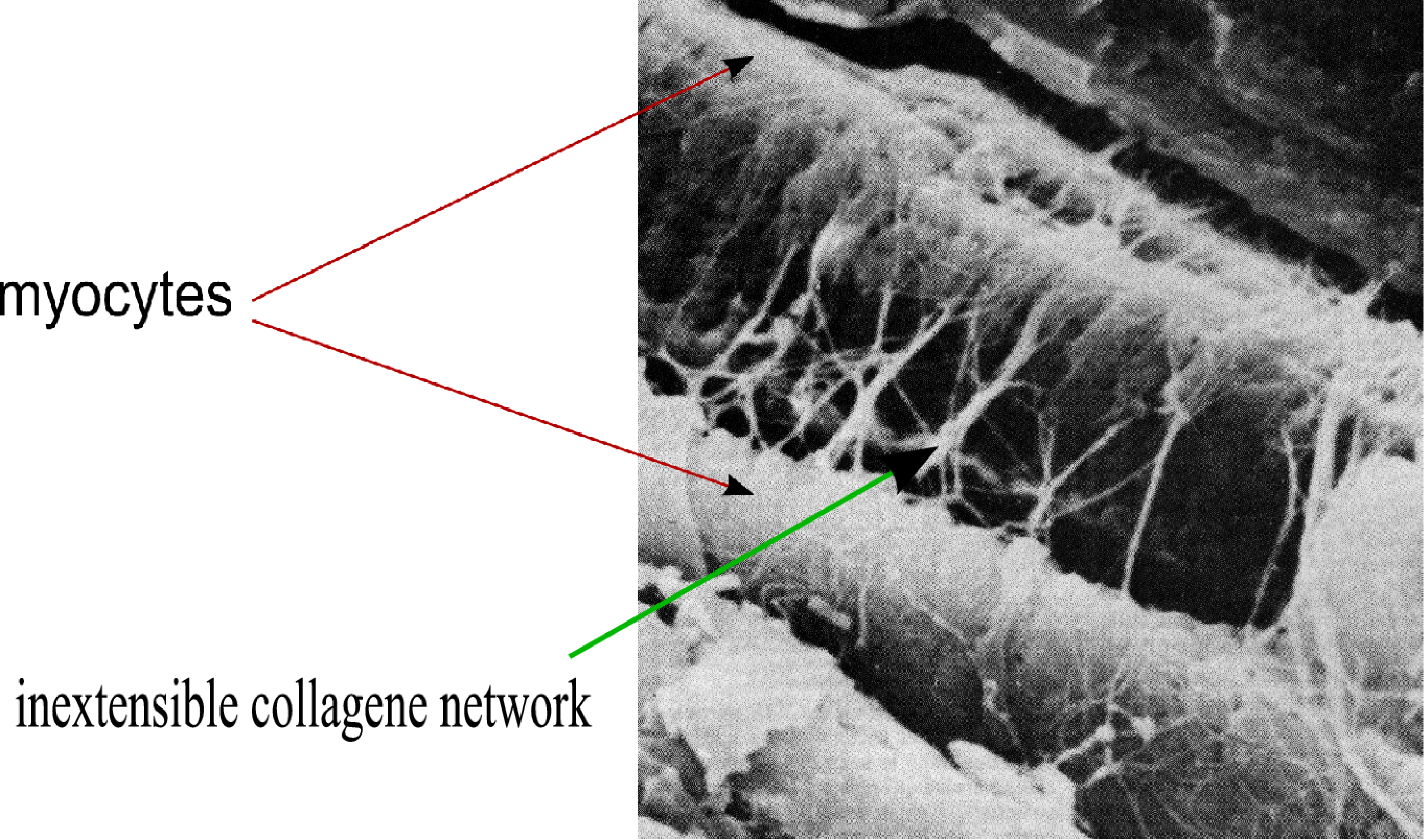}
\caption{Collagen network surrounding the myocytes.}\label{col}
\end{figure}


\section{Constitutive law in continuum mechanics}

\subsection{Coupling between muscle fibers and collagen network}

Extrapolations from muscle fiber arrangement to myocardial stress are realistic when also
taking account the effect of the connective tissue. We believe that a
part of that connective tissue, surrounding group of myocytes, is responsible for active tension
developed in the perpendicular direction of the muscle fibers running on the tangential plane of the
ventricular wall.

Based on the previous SEM observations, we proposed a connective tissue
organization illustrated on Figure 1. We assumed that the myocytes are roughly cylindrical and that
groups of myocytes are surrounded by inextensible collagen networks. So, during the
contraction, the myocytes diameter increases and because the collagen network is inextensible,
the adjacent muscle cells become closer. Thus the pseudo-active kinematic relation between
the muscle fiber and cross-fiber extension ratios (noted $\l_f$ and $\l_{cf}$, respectively) is
$h(\l_f, \l_{cf}) = 0$ with:
\begin{equation}\label{kc}
h(\l_f, \l_{cf}) =  1 - \l_{cf}  + (\pi-2)(1-\l_f^{ -1/2})\frac{a}{D}
\end{equation}

with $D=4a+d$ where $a$ is the initial myocyte radius and $d$ is the distance between the two cells.

\begin{figure}[H]
\centering\includegraphics[scale=0.4]{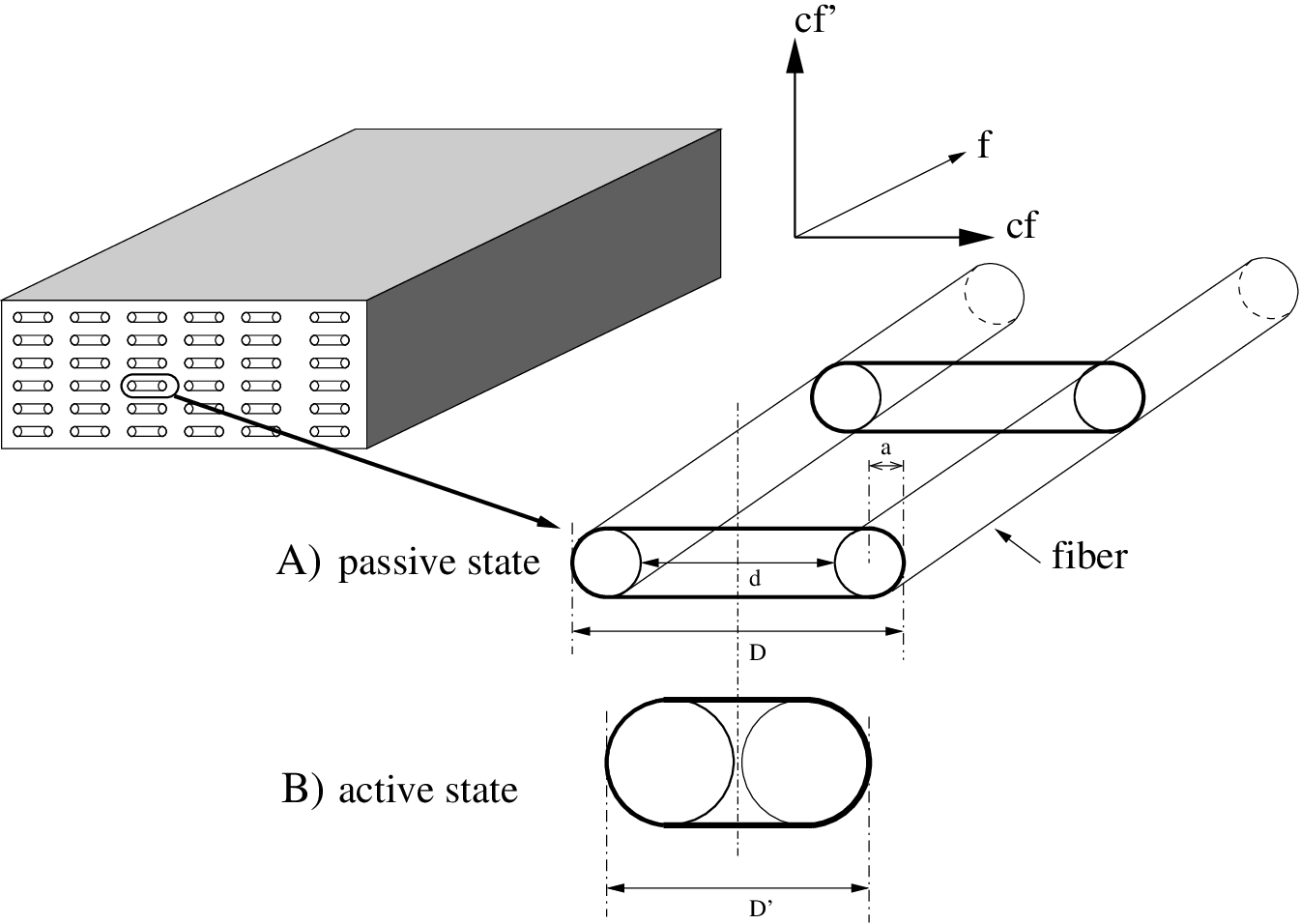}
\caption{Schematic illustration of the internal pseudo-active kinematic constraint induced
by the collagen network surrounding the myocytes. A) Before contraction. B) After or during
contraction.}\label{fib}
\end{figure}

\subsection{Constitutive law for the myocardium under internal
pseudo-active kinematic constraint}

\begin{figure}[H]
\centering\includegraphics[scale=0.45]{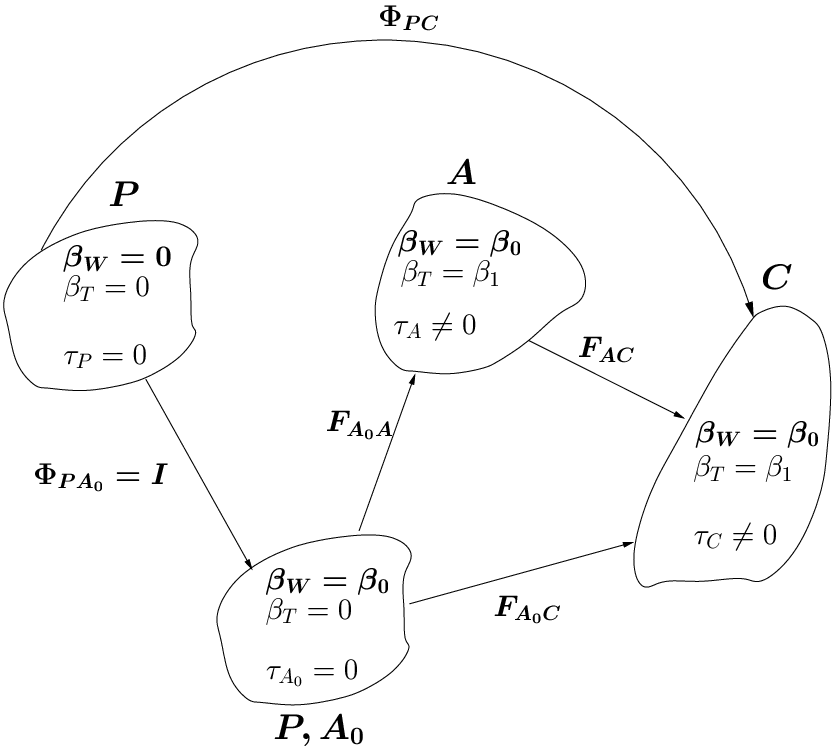}
\caption{Description of the active rheology approach.}\label{rheo}
\end{figure}

To be consistent with our mathematical formulation, the letter $\mathbf{\Phi}$  is used for non
elastic gradient tensor and the letter $\mathbf{F}$ is used for elastic gradient tensor. The activation of the
muscle fibers changes the properties of the material and at the same time contracts the muscle
itself. To have a continuous elastic description during the activation of the myocardium, we used an
approach similar to the one proposed by Chadwick \cite{C82}, Ohayon and Chadwick \cite{OC88}, 
Taber \cite{T91}. From its passive zero-stress state $P$, the free
activation of the muscle fibers is modeled by the following two transformations (Fig.2): the first
one (from state $P$ to virtual state $A_0$) changes the material properties without changing the
geometry, and the second one (from $A_0$ to $A$) contracts the muscle without changing the
properties of the material. Thus, the former is not an elastic deformation and is described by the
gradient tensor $\mathbf{\Phi}_{PA_0}=\mathbf{I}$ where $\mathbf{I}$ is the identity matrix. In that
first transformation, only the strain energy function is modified using an activation function
$\beta$, which may depend on the cardiac cycle time and some ionic concentration (calcium for instance). 
The second transformation is an elastic deformation caused only by the active tension delivered by the 
fibers and takes care of  the internal
kinematic constraint (Eq.(\ref{kc})). This last transformation is described by the gradient tensor
$\mathbf{F}_{A_0A}$. Thus the transformation from state $P$ to state $A$ is a non elastic
transformation  ($\mathbf{\Phi} _{PA}  =\mathbf{\Phi} _{PA_0 } \mathbf{F}_{A_0 A}  $), but can be
treated mathematically as an elastic one because  $\mathbf{\Phi} _{PA}  = \mathbf{F}_{A_0 A}$.
Finally, external loads are applied to state $A$ deforming the body through into $C$ (Fig. \ref{rheo}).

The change of the material properties of the myocardium during the cardiac cycle is described by a
parameter-dependent strain-energy function per unit volume of state $P$ noted $W(\mathbf{E}_{PH})$:

\begin{equation}\label{SEF}
W(\mathbf{E}_{PH} )=-\frac{1}{2}\,p_H(I_3(\mathbf{E}_{PH})-1)+W^*(\E_{PH} ) +
\delta_{AH}\,W_{\stackrel{pseudo}{ active}}(\E_{PH}) \end{equation}

with
\begin{equation}\label{SEF*}
W^*(\mathbf{E}_{PH} )=W_{pas}(\E_{PH} ) + \beta \,W_{act}^f(\E_{PH})
\end{equation}

where $\E_{PH}$ is the Green's strain tensor at an arbitrary state $H$ calculated from the zero
strain state $P$ (the state $H$ could be one of the states $A_0$, $A$ or $C$ shown in figure
\ref{rheo}), $ p_H$ is the Lagrangian multiplier resulting of the incompressibility constraint $\det
\bPhi_{PH}=1$ (see \cite{F65, M69,GZ92}),  $I_3(\mathbf{E}_{PH})$ is the determinant of the
right Cauchy-Green strain  tensor $\bC_{PH}$ ($\bC_{PH}=2\E_{PH}+\I$),   $W_{pas}$ represents the
contribution of the surrounding collagen matrix and of the passive fiber components, $W_{act}^f$
arise from the active component of the embedded muscle fibers, and   $ \beta$ is an activation
parameter equal to zero at end-diastolic state and equal to one at end-systolic state
($0\leq\b\leq 1$). The scalar $\delta_{AH}$ is equal to one if state $H$ is the state $A$ and
zero if the two states $H$ and $A$ are distinct.  The term $W_{act}^f(\E_{PH})$ gives the variation
of the muscle fibers properties during the cardiac cycle. The pseudo-active strain energy function
expressed in the last term of the right hand side of the Eq.(\ref{SEF}) is introduced in order
to satisfy the kinematic condition (Eq.(\ref{kc})) and is given by:

\begin{equation}\label{pa}
W_{\stackrel{pseudo}{ active}}(\E_{PH})=-\frac{1}{2}q_H\,h(\E_{PH})
\end{equation}

 The scalar $q_H$ introduced in Eq.(\ref{pa}) serves as an additional indeterminate Lagrange
multiplier which contributes to the pseudo-active stresses at state $H$ in fiber and the cross-fiber
directions, and $h(\E_{PH})$ is the function defined in Eq.(\ref{kc}), which may be rewritten as:

\begin{equation}\label{kcbis}
\,h(\E_{PH})=1-I_6^{1/2}+(\pi-2)\,(1-I_4^{-1/4})\,\frac{a}{D}
\end{equation}

where $I_4$ and $I_6$ are two strain invariants given by  $I_4(\E_{PH})=\f_P\cdot\bC_{PH}\cdot\f_P$
and  $I_6(\E_{PH})=\f^\bot_P\cdot\bC_{PH}\cdot\f^\bot_P$  in which the fiber and the  perpendicular
fiber directions (this last one corresponding to the direction of the collagen struts) are respectively characterized in  state $P$ by the unit vectors  $\f_P$ and
$\f^\bot_P$. In an arbitrary  deformed state $H$, the direction of these two unit vectors  are noted
$\f_H$ and $\f'_H$ and are respectively defined by:
$$\f_H=\frac{\mathbf{\Phi}_{PH}\cdot\f_P}{\norm{\mathbf{\Phi} _{PH}\cdot\f_P}} \;\mbox{ and } \; 
\f'_H=\frac{\mathbf{\Phi} _{PH}\cdot\f^\bot_P}{\norm{\mathbf{\Phi} _{PH}\cdot\f^\bot_P}} \quad .$$
The tensor $\bC_{PH}$ is the right
Cauchy-Green strain  tensor ($\bC_{PH}=2\E_{PH}+\I =\mathbf{\Phi}^T _{PH}\mathbf{\Phi} _{PH}$ ). The
superscript `$T$' is used for the transpose matrix and $\Vert\cdot\Vert$ stands for the euclidian norm.
Note that $I_4$ and $I_6$  are directly related respectively to the fiber and cross-fiber extension
ratios (we have $I_4=\l_f^2$ and $I_6=\l_{cf}^2$). In our notations $\l_f$ is related to the fiber
direction $\f_H$ and $\l_{cf}$ to the cross-fiber direction $\f'_H$ (Figure \ref{fib}). We treat the
myocardium as a homogeneous, incompressible, and hyperelastic material transversely isotropic with
respect to the local muscle fiber direction. 

In this study, the passive strain-energy function
is \cite{LY98}
\begin{equation}\label{SEFLY}
  W_{pas}(\E_{PH})=C_1^p (e^Q  - 1)
\end{equation}
\begin{equation}\label{Q}
{\rm with }\qquad Q=C_2^p (I_1  - 3)^2  + C_3^p ( I_1  - 3)(I_4  - 1) + C_4^p (I_4
- 1)^2
\end{equation}
For the active strain-energy we modified the function found by Lin and Yin \cite{LY98} by substracting
the ``beating term" $C_5^a$:
\begin{equation}\label{Wact}
W_{act} (E_{PH})=C_1^a(I_1  - 3)(I_4  - 1) + C_2^a (I_1  - 3)^2
+ C_3^a (I_4  - 1)^2  + C_4^a (I_1  - 3)
\end{equation}
where $C_i^p,\quad i=1,\cdots,4$ and $C_i^a,\quad i=1, \cdots 4$ are material constants and
$I_1$ is the first principal strain invariant given by  $I_1(\E_{PH})=tr\,\bC_{PH}\;$. 

The beating term is defined as the part of the active strain-energy function responsible for the change of
geometry when the muscle is activated and submited to no external loading. To incorporate the beating
behavior, the parameter-dependent beating tension $\beta\,T^{(0)}$ was applied in the deformed
fiber direction. In our approach, the active loaded state $C$ of the myocardial tissue is
obtained in two steps. In the first step and at a given degree of activation $\beta$,  we derived
and quantified the internal pseudo-active stresses by looking the free contraction configuration of
the tissue (state $A$, Figure  \ref{rheo}). Then, in a second step we applied the loads on the
active myocardial tissue under the internal pseudo-active stresses previously found.\\[2mm]
\textit{\textbf {Step 1: determination of the free contraction state $A$-}}
During the cardiac cycle and  at a given degree of activation $\beta$, the Cauchy stress tensor in
state $A$ (noted  $\bt_A$) is given by:

\begin{equation}\label{CauchyA}
\bt_A   = - p_A \I+ \bPhi _{PA} \frac{{\partial W^*\left( {\E_{PA}} \right)}}{{\partial \E_{PA}
}}\bPhi _{PA}^T  + \beta T^{( 0)} \f_A  \otimes \f_A+ \bt_A^{\stackrel{pseudo}{ active}}
\end{equation}
\begin{equation}\label{taupa}
\hbox{with}\qquad\bt_A^{\stackrel{pseudo}{ active}} = \bPhi _{PA} \frac{{\partial
W_{\stackrel{pseudo}{ active}}\left( {\E_{PA}} \right)}}{{\partial \E_{PA} }}\bPhi _{PA}^T
\end{equation}
where  the symbol $\otimes$ denotes the tensor product. The postulated mechanical coupling law
(Eq.(\ref{kcbis})) induces, during the contraction, a pseudo-active stress tensor:
\begin{equation}\label{taupseudo}
 \bt_A^{\stackrel{pseudo}{ active}}= T_A^f\, \f_A  \otimes \f_A+ T_A^{cf}\,\f'_A  \otimes \f'_A
 \end{equation}
These two stress tensor components $T_A^f$ and $T_A^{cf}$ are activation-dependent and behave as
some internal tensions  in the fiber and cross-fiber directions of unit vectors $\f_A$ and $\f'_A$,
respectively. These  pseudo-active tensions  are defined by:
\begin{equation}\label{T}
T_A^f=2\, \frac{{\partial W_{\stackrel{pseudo}{active}}}}{{\partial
I_4 }}\,\norm{\bPhi _{PA}\cdot\f_P}^2 \quad \hbox{;} \quad
T_A^{cf}=2\, \frac{{\partial W_{\stackrel{pseudo}{ active}}}}{{\partial
I_6}}\,\norm{\bPhi _{PA}\cdot\f^\bot_P}^2
\end{equation}\\[2mm]
\textit{\textbf {Step 2: determination of the physiological active loaded state $C$-}}
These previously found internal pseudo-active tensions $ T_A^f$ and $T_A^{cf}$  were introduced
in the expression of the stress tensor at loaded state $C$. Therefore,  at a given degree of activation $\beta$, 
the Cauchy stress tensor in the physiological state $C$ (noted $\bt_C$) is given by:
\begin{equation}\label{CauchyC}
\bt_C  = - p_C \I+ \bPhi _{PC} \frac{{\partial W^*\left( {\E_{PC}} \right)}}{{\partial \E_{PC}
}}\bPhi _{PC}^T  + \left( \beta T^{( 0)}+T_A^f\right) \, \f_C\otimes \f_C+
T_A^{cf} \, \f'_C\otimes \f'_C
\end{equation}
The  suggested constitutive law for the active myocardium (Eqs.(\ref{SEF})-(\ref{CauchyC}))
allows to simulate  the left ventricle behavior during the whole cardiac cycle.
Thus, in this  law: (i)  the anisotropic behavior is incorporated in the expressions of passive,
active and pseudo-active strain energy functions by the terms  $I_4$ and $I_6$, (ii) the kinematic
contraction is accounted for by a beating tension  $\beta\,T^{(0)}$ in the fiber direction, (iii)
the change of properties is expressed by the active  strain  energy term $\beta \,W_{act}$,
and (iv) the coupling effect between the collagen network and the MF is accounted for
by the two internal pseudo-active tensions $ T_A^f$ and $T_A^{cf}$ in the fiber and cross fiber
directions $\f_C$  and $\f'_C$, respectively.


\section{Variational formulation and finite element method}


 The undeformed body state $P$ contains a volume $V$ bounded by a closed surface $\mathcal{A}$, and
the arbitrary deformed body state is, as before, noted $H$. The corresponding position vectors, in
cartesian base unit vectors, are $\R = Y^R \e_R$  and $\r = y^r \e_r$ , respectively. However, we
write the equations with suitable curvilinear systems of world coordinates noted $\Theta^{\rm A}$ in
the reference configuration (state $P$) and $\theta^\a$ in the deformed configuration (state $H$): see Fig. \ref{coord} in appendix \ref{COORD}.
In this paper we use the same conventional notations (Table 1 in appendix \ref{COORD}) for
vectors, tensors and coordinates systems  than Costa et al. \cite{CHWWGM961, CHWWGM962},
where:

\begin{itemize}
  \item[-] Capital letters are used for coordinates and indices of tensor components associated to
state P, and lower case letters are related to state $H$.
  \item[-] $\mathbf{G}$ and $\mathbf{g}$ are the base vectors in states $P$ and $H$,
respectively, for which parenthetical superscript indicates the
associated coordinate system (for example $\G_I^{(x)}=
\ds{\partial \R}/{\partial X^I}=\R_{,I}^{(x)}$  and
$\mathbf{g}_i^{(x)}=\ds{\partial \r}/{\partial x^i} =
\r_{,i}^{(x)}$).
\end{itemize}

The Lagrangian formulation of the virtual works principle is given
by (\cite{CHWWGM961, M69})
\begin{equation}\label{LF}
\int\limits_V {P_H^{IJ} \Phi _J^{\cdot\alpha}\nabla _I}(\delta
u_\alpha)\, dV = \int\limits_V {\rho (b^\alpha-\gamma
^\alpha)}\delta u_\alpha \, dV +(1-\delta_{AH})\, \int\limits_{A_2 }
\s.\bd\mathbf{u}\, dA
\end{equation}

where $P_H^{IJ}$ are the components of the second Piola-Kirchhoff stress tensor at state $H$,
$\P_H$, referred to the base tensor $\G_I^{(x)}\otimes\G_J^{(x)}$, $\Phi
_I^{\cdot\alpha}=\ds{\partial \theta^\alpha/\partial X^I}$  are the components of the gradient
tensor $\bPhi_{PH}$  in the base tensor $ \mathbf{g}_\alpha ^{(\theta )}\otimes\G^{(x)I} $,
$\bd\mathbf{u}=\delta u_\alpha \mathbf{g}^{(\theta )\alpha }$  is an arbitrary  admissible
displacement vector, $\nabla _I (\delta u_\alpha  ) = \ds{\partial \delta u_\alpha/\partial
X^I}-\mathbf{g}_{\alpha ,I}^{(\theta )} \cdot\mathbf{g}_{}^{(\theta )\beta } \delta u_\beta$ are the
components of the covariant differentiation vector $\bd\mathbf{u}$ in the base vectors
$\mathbf{g}^{(\theta )\alpha}$ (i.e. $\nabla _I (\delta u) = \nabla _I (\delta u_\alpha
)\mathbf{g}_{}^{(\theta )\alpha }$). The previous differentiation is done with respect to the
locally orthonormal body coordinates ( $X^I,\; I=1,2,3$) which coincide with the local muscle fiber
direction in state $P$. The material density in the undeformed body state $P$ is $\rho$,
$\mathbf{b}=b^\alpha\mathbf{g}_\alpha ^{(\theta )}$ is the body force vector per unit mass, $\bg  =
\gamma ^\alpha  \mathbf{g}_\alpha ^{(\theta )}$ is the acceleration vector, $\s$ is the surface
traction per unit area of $\mathcal{A}$, and $A_2$ is the part of $\mathcal{A}$ not subject to
displacement boundary conditions.
The Lagrangian formulation for incompressibility is given by
\begin{equation}\label{incomp}
\int_V\left(\det g^{(x)}_{IJ}-1\right)\,p^*\,dV=0
\end{equation}
where the metric tensor $g^{(x)}_{IJ}$ is  defined in table 1,
and $p^*$ is an arbitrary admissible pressure.
Lastly   the  Lagrangian formulation for the additional pseudo-active kinematic constraint is given
by
\begin{equation}\label{kcform}
\delta_{AH}\,\int_V h(I_4,I_6)\,q^*\,dV=0
\end{equation}
for all admissible $q^*$. Eqs.(\ref{LF})-(\ref{incomp}) -(\ref{kcform}) represent the variational formulation of a system of
nonlinear partial differential equations.
For an incompressible medium ($\det\bPhi_{PH}=1$), the relation between the second Piola-Kirchoff
stress tensor $\P_H$ and the Cauchy stress tensor $\bt_H$  is (\cite{M69})
\begin{equation}\label{Piola}
    \P_H=\bPhi_{PH}^{-1}\,.\,\bt_H\,.\,(\bPhi_{PH}^{-1})^T
\end{equation}
A complete expression of the components of   $\P_H$ in both states $H=A$ and $H=C$ are given
in appendix \ref{P}. The surface  traction per unit of undeformed area of $\mathcal{A}$,
$\mathbf{s}=s^\a\mathbf{g}_\a^{(\theta)}$,  is a known loading boundary which could be written
using physical Cauchy stress.


\subsection{Finite element approximation}

\begin{figure}[H]\label{mesh}
\centering\includegraphics[scale=0.6]{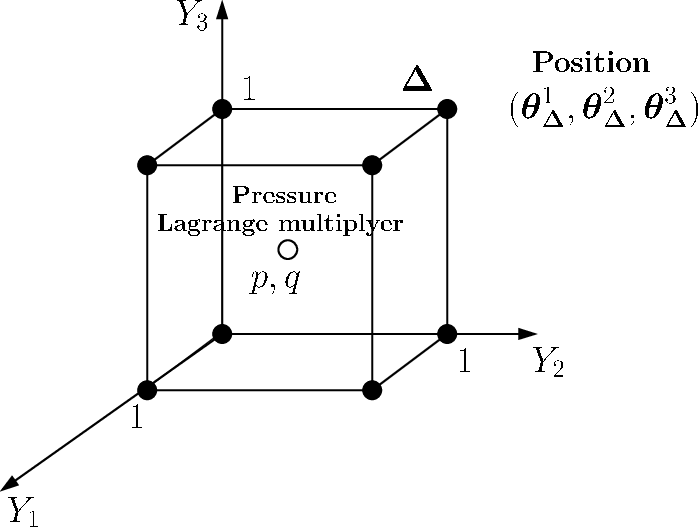}
\caption{$Q1-Q0$ element for displacements, pressure and pseudo-active Lagrange multiplyer.}
\end{figure}

Through this paper, we use a three dimensional finite element with Lagrange trilinear interpolation
for the displacements and uniform Lagrangian multipliers to compute an approximate solution of
Eqs.(\ref{LF})-(\ref{incomp})-(\ref{kcform})  on a rectangular mesh (see Fig. \ref{mesh}), where we neglect the
acceleration and body forces ($\mathbf{b}=0$, $\bg=0$). This element is commonly used
and is relevant for the finite element approximation of this type of problem where 
kinematics constraints must be satisfied (for more details see \cite{O72, GT841, GT842, GT89,C80, QA94}). 

Let $(\xi_K)$ the Lagrangian normalized finite element coordinates (Figure
\ref{coord}),  the deformed geometric coordinates $\theta^\a$ in element $e$ are interpolated as

\begin{equation}\label{inter}
\theta^\a=\sum_{n(e)=1}^{8}\psi_{n(e)}(\xi_1,\xi_2,\xi_3)\,
\theta^\a_{n(e)}
\end{equation}

where $\psi_{n(e)}$ is the base function associated with the local node $n(e)$ and
$\theta^\a_{n(e)}$ is the $\a$-coordinate of the local node $n$ of element $e$.\\[12pt]
Let $\O_\Delta^{n(e)}$ be the connectivity matrix defined by
\begin{equation}\label{connect}
\O_\Delta^{n(e)}=  \left\lbrace
\begin{array}{rcl}
  1 & \hbox{ if } & \Delta(n(e),e)=\Delta \\
  0 &  & \hbox{otherwise}
\end{array}
\right.
\end{equation}

The FE approximation of Eqs.(\ref{LF})-(\ref{incomp})-(\ref{kcform})  is

\begin{eqnarray}
\sum_e \sum_{n(e)=1}^{8}\O_\Delta^{n(e)}\int_{V_e} P_H^{IJ}\,
\Phi_J^{\cdot\a}\,\nabla _I(\psi_{n(e)})\, dV &=&(1-\delta_{AH})\,\sum_e
\sum_{n(e)=1}^{8}\O_\Delta^{n(e)}\int_{A_{2_e}}
s^\a\,\psi_{n(e)}\, dA    \label{FEA}\\[12pt]
\forall e \,,\, \int_{V_e}  \left(\det g^{(x)}_{IJ}-1\right)\,dV&=&0
\label{incomp2}\\[12pt]
\forall e \,,\, \delta_{AH}\,\int_{V_e}h(I_4,I_6)\,dV=0\label{kcform2}
\end{eqnarray}

with $\Delta=1,\cdots,\Delta_{max}$, $\a=1,2,3$, where $V_e$ is the volume of the element $e$, $A_{2_e}$ is the part of  $A_e$ (boundary of the element $e$) non subject to displacement  conditions.%
%
\subsection{Finite element solution method}
We proceed in two steps. The first one consists in the determination of the pseudo- active stresses
$T_A^f$ and $T_A^{cf}$ as  functions  of the activation parameter $\beta\in [0,1]$ by looking for
the state $A$ ($\delta_{AH}=1$). We solve the system (\ref{FEA})-(\ref{incomp2})-(\ref{kcform2})
with zero right hand side for (\ref{FEA}) (free active contraction) and $P_A^{IJ}$ given by
Eq.(\ref{PIJA}) . The unknowns of this nonlinear system of equations  are $\left(\theta^\alpha_{\Delta},
p_A(e), q_A(e)\right)$  with $\alpha=1,2,3$, $\Delta=1,\cdots,\Delta_{max}$ and  $ e=1,\ldots,
e_{max}$ where $e_{max}$ is the total number of elements involved in the mesh.  We derive
$T_A^f$ and $T_A^{cf}$ for a given $\beta$ according to Eq.(\ref{T}). Then in a next step  we can
compute any physiological active loaded state $C$ solving the system (\ref{FEA})-(\ref{incomp2})
with  $\delta_{AH}=0$ and $P_C^{IJ}$ given by Eq.(\ref{PIJC}). To solve the sytem in both cases we use the Powell method \cite{P70} implemented in the package \textsc{minpack} \cite{MGH80}.
\section{Results and discussion}
This section is devoted to the numerical simulation of two types of material:
\begin{itemize}
 \item a thin sample of living myocardium for which a cartesian coordinate is used,
\item  an active thick-walled cylinder for which cylindrical coordinate is used.
\end{itemize}
For these two simple configurations the exact displacements are solutions of a nonlinear system and can be computed with a high degree of accuracy.
A very good agreement between the exact and the computed solutions of the finite element nonlinear system (\ref{FEA})-(\ref{incomp2})-(\ref{kcform2}) is obtained. More precisely in all the cases the $L^2$ norm of the error
is less than $10^{-09}$.

\subsection{Case of a thin sample of living myocardium}

\begin{figure}[H]
\centering\includegraphics[scale=0.35, angle=-90]{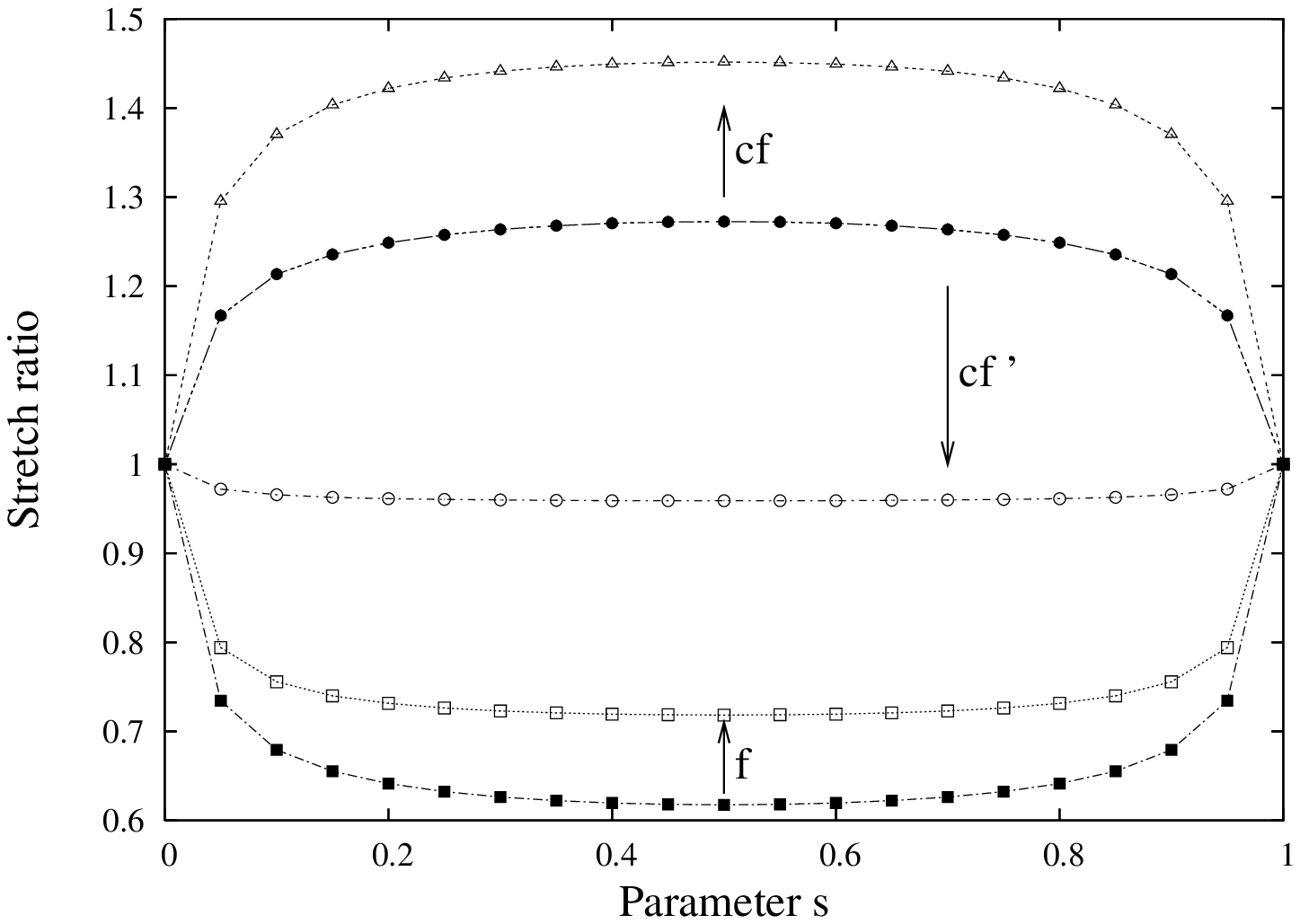}
\caption{Free contraction test with $\beta(s)=sin^2(\pi s)$: effect of the pseudo-active kinematic constraint.
The empty and full symbols indicate that the coupling effect is acting or not, respectively. The fiber and cross-fiber
directions are noted (f), (cf) and (cf~') and are defined in figure \ref{fib}. Arrows show the curve modification when the
pseudo-active kinematic constraint behaves.}
\label{free} \end{figure}

We simulated the loading of a thin sample of living myocardium ($1.0 \times 1.0 \times 0.1\; {\rm cm}^3$)
in  which the MF are uniformly oriented in one direction. The coefficients involved in the strain energy-function
are those of Lin and Yin  \cite{LY98}: $C_1^p$=0.292 kPa, $C_2^p$=0.321, $C_3^p$=-0.260, $C_4^p$=0.201,
$C_1^a$=-3.870 kPa, $C_2^a$=4.830 kPa, $C_3^a$=2.512 kPa and $C_4^a$=0.951 kPa.
For the beating tension, a good agreement between the previous experimental results and our theoretical
solution is obtained for  $T^{(0)}$=0.6 kPa. Nevertheless, the control simulation was performed with $a/D$=0.2 and
$T^{(0)}$=35 kPa. This higher value of $T^{(0)}$ is more adapted to the description of the left
ventricular performance \cite{OCJUA01}.\\[2mm]
\textit{\textbf{Influence of the collagen network on the systolic wall thickening-}}
The free contraction test is performed with no external displacement or force on the
boundaries of the sample, but just in activating the tissue. In this simulation we
used the following activation function : $\beta(s)=sin^2(\pi s)$. Compare to the case where the kinematic constraint
is not taken into account, one can see an increase of the cross-fiber extension ratio which
is in the tangential plane of the ventricular wall (Figure \ref{free}).
At the end-systolic state (i.e. when $\beta=1$), this ratio goes from the value 1.25 if we neglect the coupling effect,
to 1.45 when considering the kinematic constraint induced by the collagen. So, the connective tissue
could account for $16$ \% of normal end-systolic wall thickness. This increase is clearly dependent of the
geometrical parameter ratio $a/D$ and the maximal beating tension $T^{(0)}$.

\begin{figure}[H]
\centering\includegraphics[scale=0.35, angle=-90]{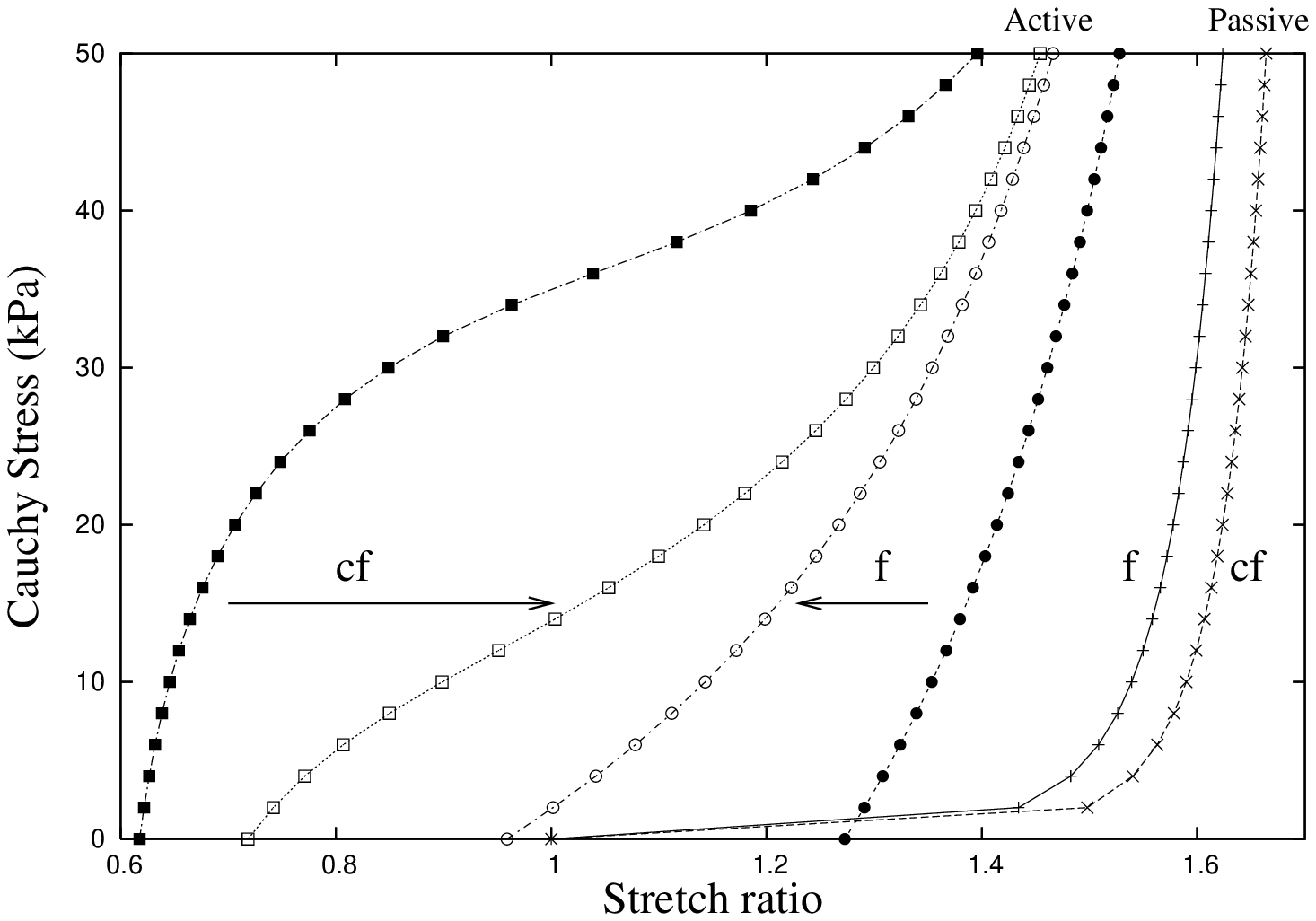}
\caption{Active and passive uniaxial extension tests: effect of the pseudo-active kinematic
constraint. The empty and full symbols indicate that the coupling effect is acting or not, respectively.
The fiber and cross-fiber
directions noted (f) and (cf) are defined in figure \ref{fib}. Arrows show the curve modification when the
pseudo-active kinematic constraint behaves.}
\label{uni} \end{figure}

\textit{\textbf{Influence of the collagen network on the pseudo-active tension-}}
Table \ref{table} shows the effect of the geometrical parameter $a/D$ and the maximal active tension $T^{(0)}$,
on the fiber and cross-fiber
stresses (noted $\sigma_{11}$ and $\sigma_{22}$, respectively). These effects were given in the case of
an equibiaxial extension loading ($\l_{f}=\l_{cf}=1.2$) of an activated sample of myocardium ($\beta=1$).
These two stresses increase with $T^{(0)}$, but are not very sensitive
to the geometrical ratio $a/D$. We can observe also, that by neglecting the interaction between the collagen
network and the MF: (i)  the cross-fiber stress is not affected by the amplitude of the beating tension, and
(ii) the stress ratio  $\sigma_{22}$/$\sigma_{11}$ decreases when $T^{(0)}$ increases. These results
mean that the usual strain-energy functions considered for the myocardium are not able to generate any transverse
pseudo-active tension.
Moreover, the results obtained for the uniaxial tests of an active or a passive sample, with or without the effect
of the collagen on the MF, are shown in Figure \ref{uni}. Because the coupling effect between the collagen and the MF
is an active mechanism, the passive stress-strain relations are not affected by the kinematic constraint.
The mechanical properties of the active tissue, in the fiber and cross-fiber directions, become comparable when the
coupling effect acts.

\begin{table}[ht]
\caption{Effect of active tension  $T^{(0)}$ and geometrical parameter $a/D$}
\label{table}
$
\begin{array}{|c|c|c|c|c|c|c|}
\hline
\hspace{3mm}
 a/d
 \raisebox{-1 ex}[1cm][0cm]{
 \hspace{33mm}
 \unitlength=1mm
 \begin{picture}(20,10)
 \line(-4,1){47}
\end{picture}}
\hspace{-40mm}
\raisebox{2.5 ex}[0cm][0cm]{$T_0 \;(KPa)$}
&& 5 & 15 & 25  & 35 & 45 \\
\hline

\hspace{20mm}&  \raisebox{2 ex}[1cm][0cm]{$\sigma_{22}/\sigma_{11} \; (\%)$}   &\raisebox{2
ex}[1cm][0cm]{  50.92} &\raisebox{2 ex}[1cm][0cm]{53.73 } &\raisebox{2 ex}[1cm][0cm]{57.65 }
&\raisebox{2 ex}[1cm][0cm]{60.25  } & \raisebox{2 ex}[1cm][0cm]{  61.90 } \\
\cline{2-7}
 \raisebox{6 ex}[0cm][0cm]{ 0.10} & \raisebox{2 ex}[1cm][0cm]{$ \sigma_{22} \; (KPa)$}  &
\raisebox{2 ex}[1cm][0cm]{5.94} &\raisebox{2 ex}[1cm][0cm]{ 11.44} & \raisebox{2
ex}[1cm][0cm]{17.78} & \raisebox{2 ex}[1cm][0cm]{ 24.32} &  \raisebox{2 ex}[1cm][0cm]{30.88} \\
\hline
\hspace{20mm}&  \raisebox{2 ex}[1cm][0cm]{$\sigma_{22}/\sigma_{11} \; (\%)$}   &\raisebox{2
ex}[1cm][0cm]{  51.44} &\raisebox{2 ex}[1cm][0cm]{  54.56} &\raisebox{2 ex}[1cm][0cm]{58.60  }
&\raisebox{2 ex}[1cm][0cm]{61.30} & \raisebox{2 ex}[1cm][0cm]{   63.01 } \\
\cline{2-7}
 \raisebox{6 ex}[0cm][0cm]{ 0.15} & \raisebox{2 ex}[1cm][0cm]{$ \sigma_{22} \; (KPa)$}  &
\raisebox{2 ex}[1cm][0cm]{5.97} &\raisebox{2 ex}[1cm][0cm]{ 11.48} & \raisebox{2
ex}[1cm][0cm]{17.80} & \raisebox{2 ex}[1cm][0cm]{ 24.31} &  \raisebox{2 ex}[1cm][0cm]{30.84} \\
\hline
\hspace{20mm}&  \raisebox{2 ex}[1cm][0cm]{$\sigma_{22}/\sigma_{11} \; (\%)$}   &\raisebox{2
ex}[1cm][0cm]{  51.94} &\raisebox{2 ex}[1cm][0cm]{  55.33} &\raisebox{2 ex}[1cm][0cm]{59.47  }
&\raisebox{2 ex}[1cm][0cm]{62.25} & \raisebox{2 ex}[1cm][0cm]{   64.04} \\
\cline{2-7}
 \raisebox{6 ex}[0cm][0cm]{ 0.20} & \raisebox{2 ex}[1cm][0cm]{$ \sigma_{22} \; (KPa)$}  &
\raisebox{2 ex}[1cm][0cm]{6.00 } &\raisebox{2 ex}[1cm][0cm]{ 11.51} & \raisebox{2
ex}[1cm][0cm]{17.79} & \raisebox{2 ex}[1cm][0cm]{ 24.25} &  \raisebox{2 ex}[1cm][0cm]{30.73} \\
\hline
\hspace{40mm}&  \raisebox{2 ex}[0cm][0cm]{$\sigma_{22}/\sigma_{11} \; (\%)$}   &\raisebox{2
ex}[1cm][0cm]{ 36.78} &\raisebox{2 ex}[1cm][0cm]{19.88  } &\raisebox{2 ex}[1cm][0cm]{13.62 }
&\raisebox{2 ex}[1cm][0cm]{10.36 } & \raisebox{2 ex}[1cm][0cm]{  8.36} \\
\cline{2-7}
 \raisebox{6 ex}[0cm][0cm]{ \mbox{No kinematic constraint}} & \raisebox{2 ex}[1cm][0cm]{$
\sigma_{22} \; (KPa)$}  & \raisebox{2 ex}[1cm][0cm]{4.32} &\raisebox{2 ex}[1cm][0cm]{4.32 } &
\raisebox{2 ex}[1cm][0cm]{4.32} & \raisebox{2 ex}[1cm][0cm]{4.32} &  \raisebox{2 ex}[1cm][0cm]{4.32}
\\
 \hline
\end{array}
$
\end{table}

\newpage

\subsection{Case of an active thick-walled cylinder}


We simulate the mechanical behaviour of an active artery  under physiological blood pressure
$P_{int}$. This artery is modelled by a thick-walled cylinder with internal radius $R_{int}=2$~mm,
external radius $R_{ext}=3.5$~mm and height $L=2$~cm. We assume that the medium is made of a
hyperelastic anisotropic material with fibers  oriented in the circumferential direction. This simulation does not take into account all the complexity of the structure of an artery, so the following results must be viewed as a first approach and we focus on some qualitative aspects, particularly the wall thickening effect.

\begin{figure}[H]
\centering
\includegraphics[scale=0.25]{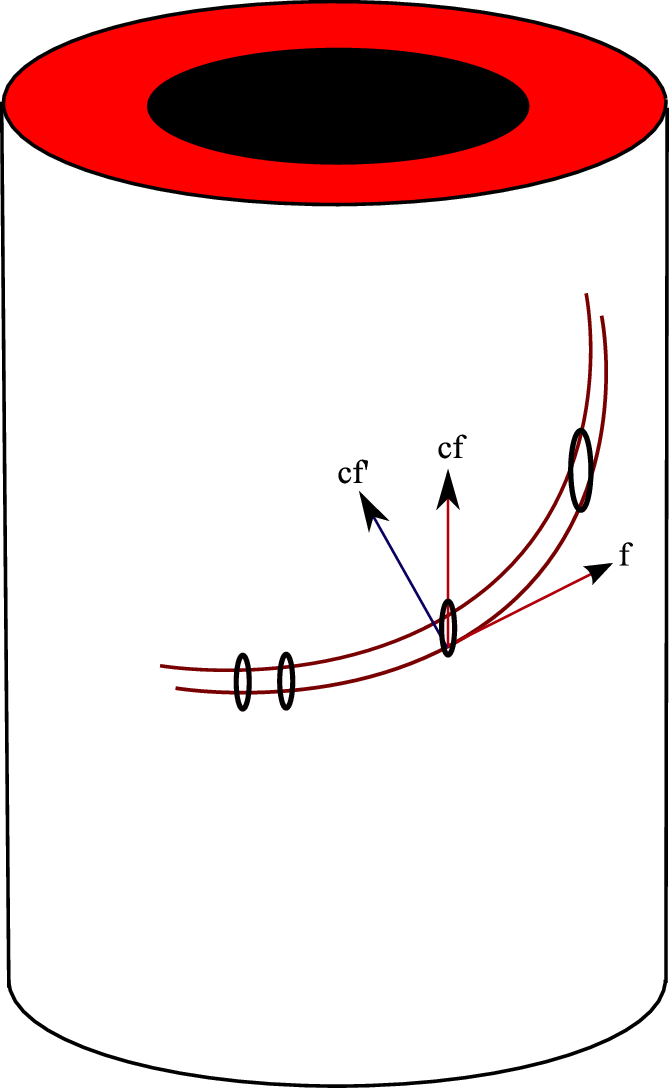}\hspace{10mm}
\includegraphics[scale=0.25]{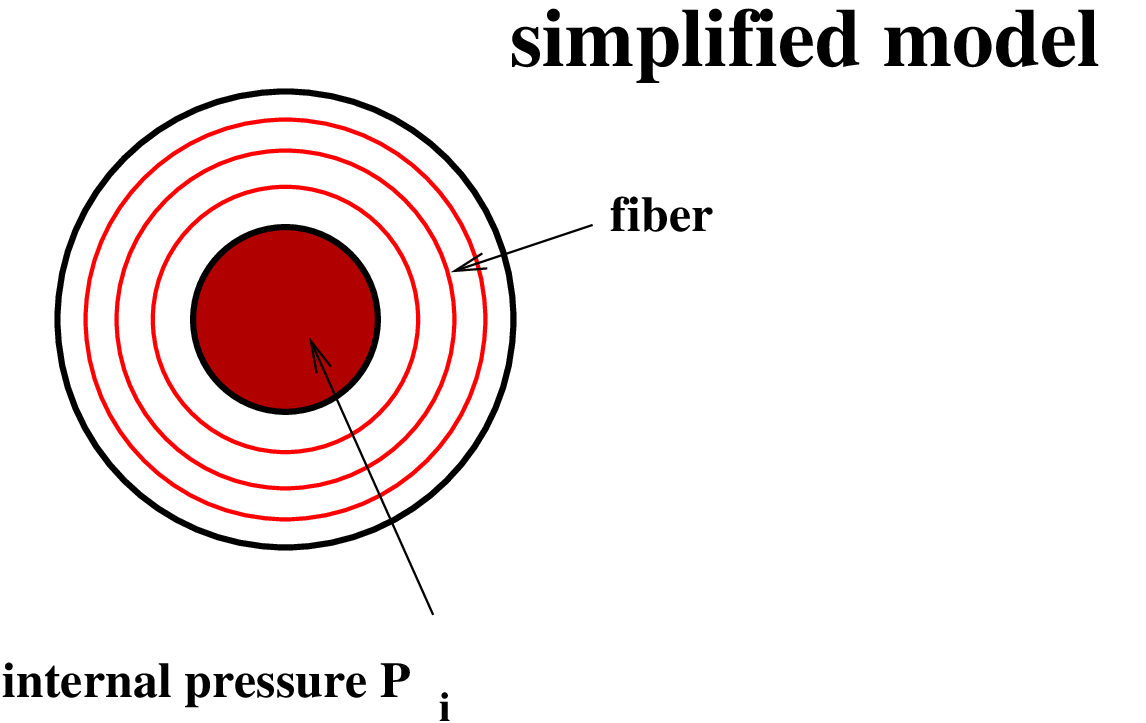}
\caption{a simplified model for the spatial organization of the  collagen fibrils. The fibrils are oriented in the $z$-direction ($cf$) and in the right part of the figure we suggest that the fibers  are oriented in the circumferential direction}
\label{artery}
\end{figure}

\begin{figure}[H]
\centering
\includegraphics[scale=0.35, angle=-90]{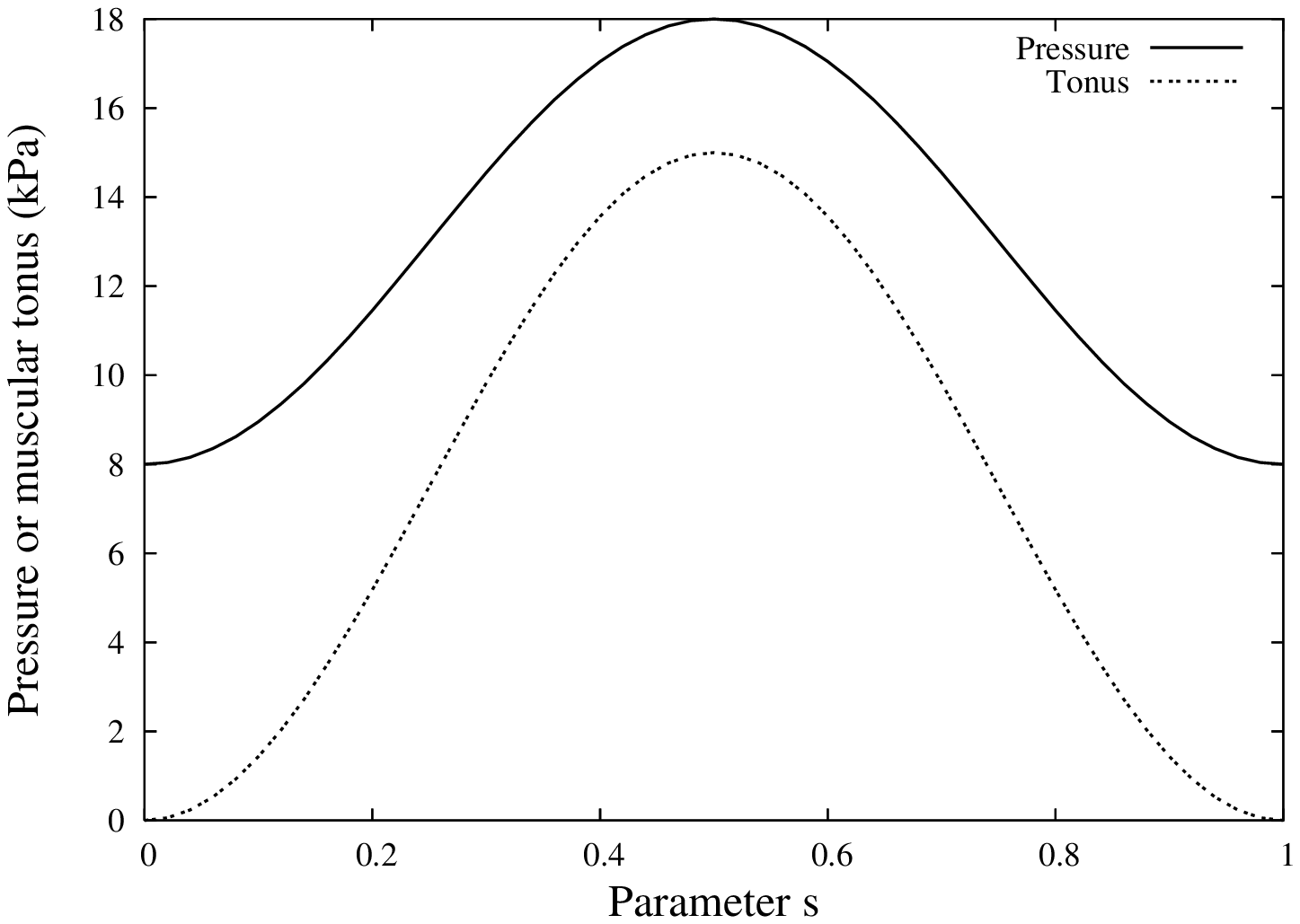}
\caption{Pressure and tonus as a function of the degree of activation}
\label{pression}
\end{figure}

\begin{figure}[H]
\centering
\includegraphics[scale=0.35, angle=-90]{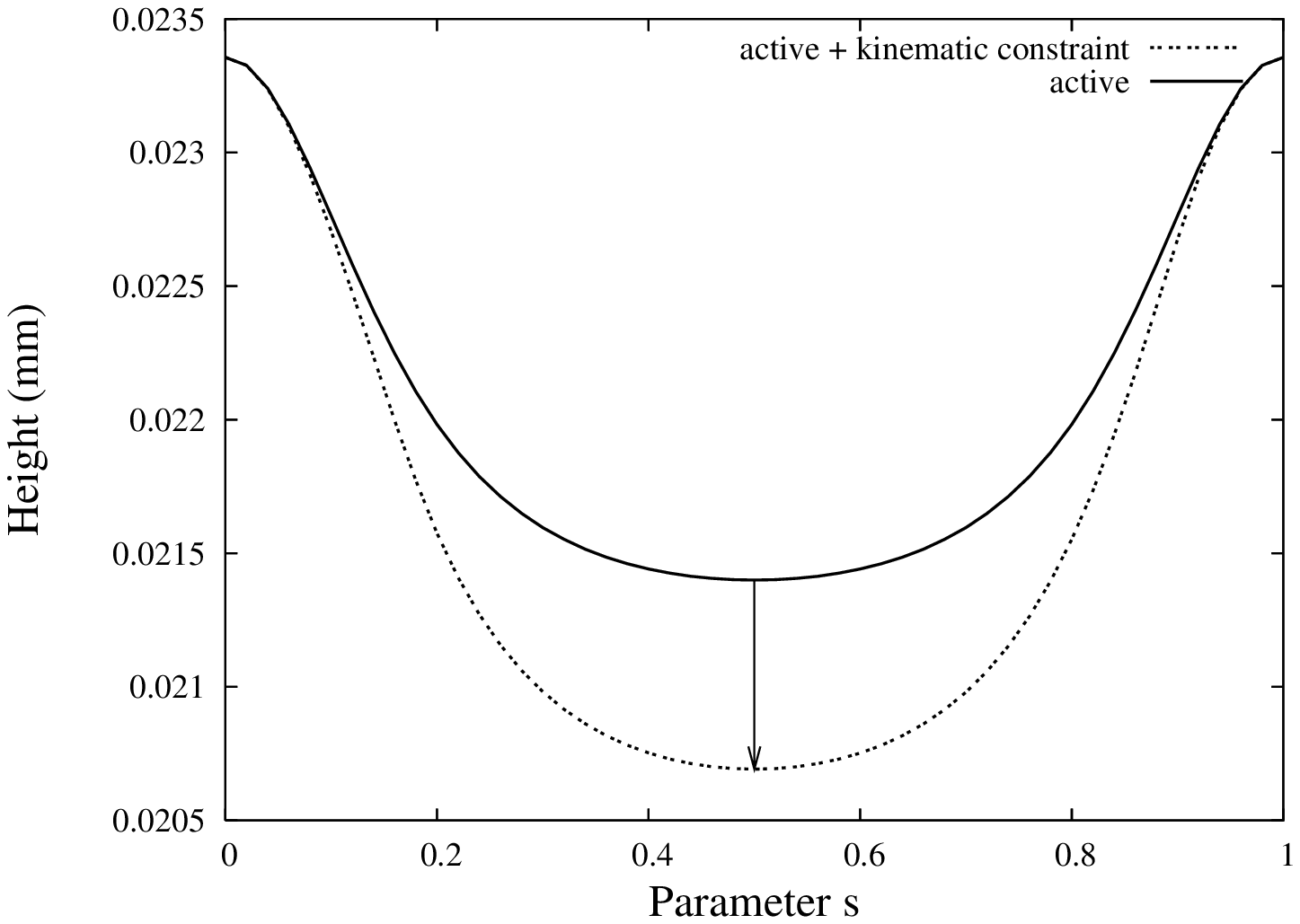}
\caption{height of the cylinder as a function of the degree of activation}
\label{hauteur}
\end{figure}

\begin{figure}[H]
\centering
\includegraphics[scale=0.35, angle=-90]{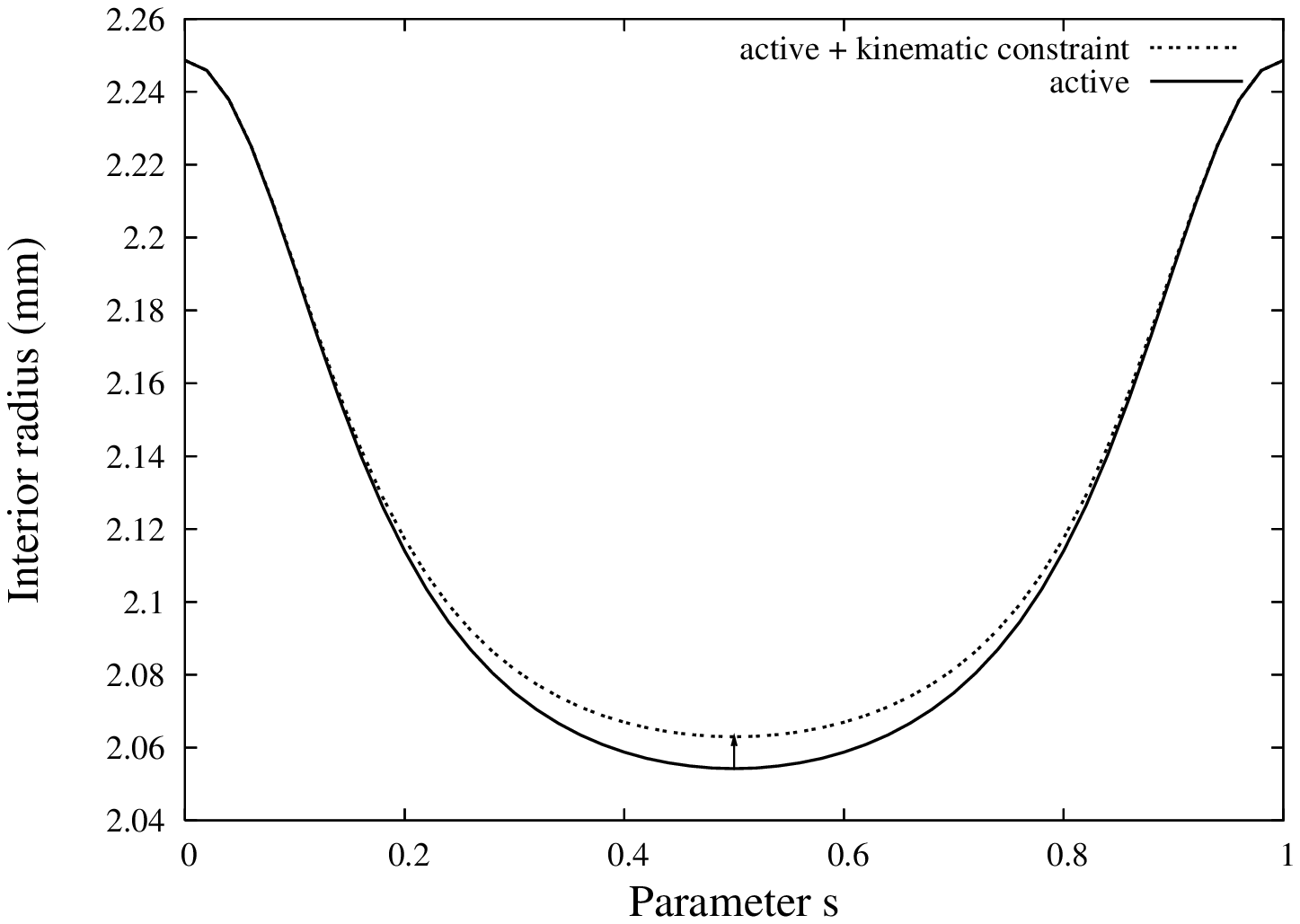}
\caption{interior radius of the cylinder as a function of the degree of activation}
\label{intrad}
\end{figure}

\begin{figure}[H]
\centering
\includegraphics[scale=0.35, angle=-90]{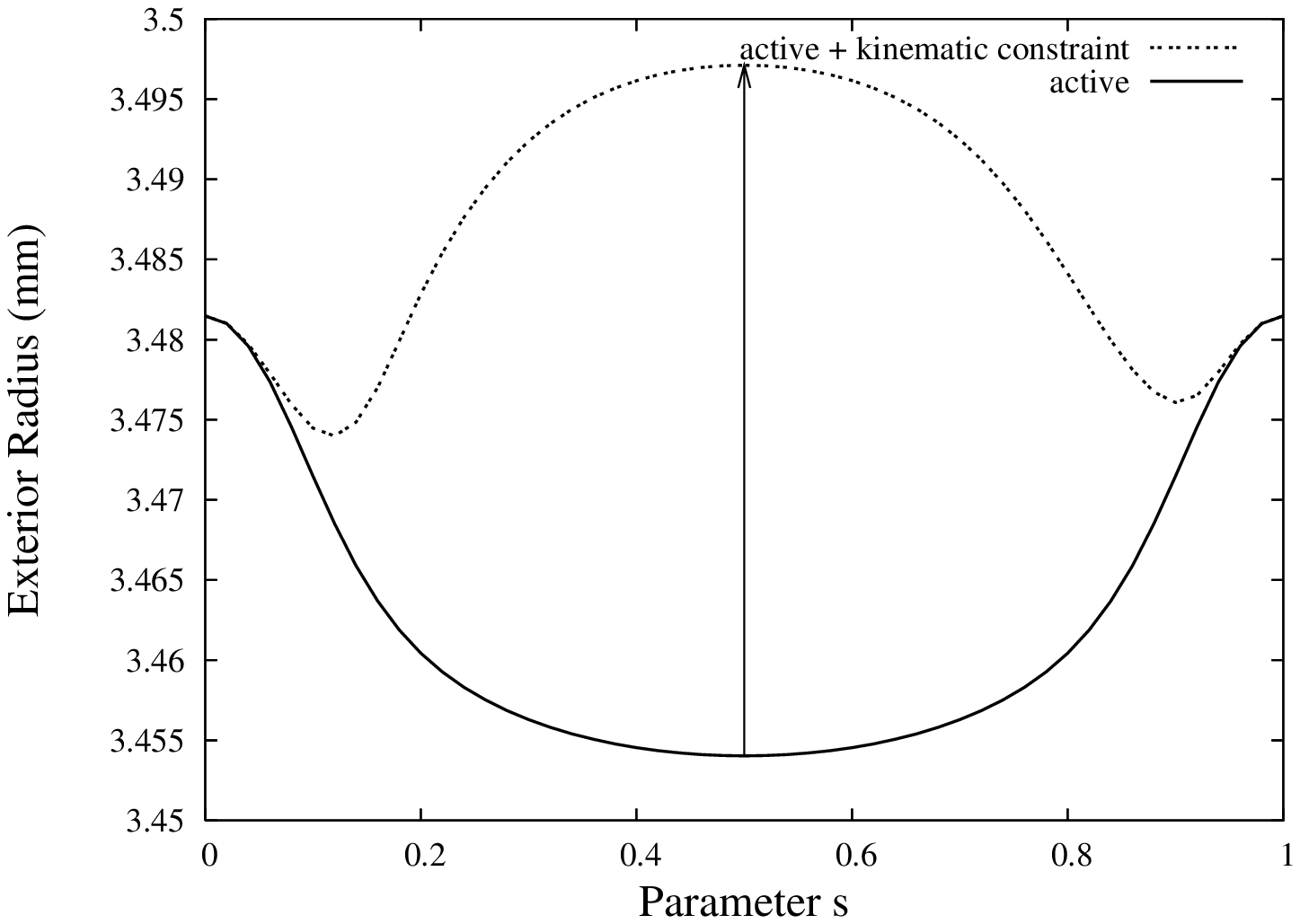}
\caption{exterior radius of the cylinder as a function of the degree of activation}
\label{extrad}
\end{figure}

In 1902, Bayliss suggested that the distension of  the vessel by blood pressure could act as a
mechanical stimulus to the vascular smooth muscle cells, thereby contributing to their tone
\cite{B02}. However, conclusive experimental support for this concept was available only  recently. We now know that the degree of vascular distension appears to be a factor of importance in determining vascular tone. We used the suggested constitutive law to model a hypothetical
autoregulation mechanism.\\
For this simulation, the active fiber tension as well as the rheological
change are in phase with the pulsatile pressure, and we use as input data the following functions:
$\b(s)=\sin^2\pi s$ with $P_{int}=8+10\,\sin^2\pi t$ (Fig. \ref{pression}).  
The resulting variations of the thick-walled
cylinder radii are presented in Figures \ref{intrad}-\ref{extrad} for this autoregulation law based on fluid
pressure. The autoregulation is defined as the relationship between the activation function $\b(s)$
and the pulsatile blood pressure $P_{int}$. Very interestingly, the results show that the kinematics
of the arterial wall may be more sensitive to the change of mechanical properties than to the blood
pressure. In other words, it appears that if there is no kinematic constraint due to the fibrils, the internal and external radii increase when the blood
pressure decreases. In fact, during this decrease of pressure, we assume that the material becomes
more compliant. Thus, the wall kinematics is mainly driven by the change of rheology. Furthermore,
although  the pressure and activation are in phase, you can create with this autoregulation law
some delay in the kinematic response. Therefore we believe that the pressure-activation interaction is a fundamental mechanism which must be well modeled to describe accurately the behaviour of the
arterial wall under physiological or pathological conditions \cite{TOH97, R93}.\\ 
On the other hand, when the kinematic constraint is activated, we observe that the interior radius still decreases but less than previously, due to the residual constraints $T_A^{f}\mbox{ and }T_A^{cf}$ (see Fig. \ref{intrad}).
The more important effect is observed for the exterior radius (Fig. \ref{extrad}),
due to a greater residual constraints $T_A^{f}\mbox{ and }T_A^{cf}$ and  the absence of exterior stress.
The global result is a wall thickening effect and a contraction in the $z$-direction due to the incompressibility condition (see Fig. \ref{hauteur}) 

\section{Conclusion}
This study shows that the connective tissue skeleton in the normal and pathological left ventricle
may have a large influence on the cardiac performance. A new constitutive law has been developed for large
deformations of an incompressible hyperelastic, and anisotropic living myocardium.
This work is based on the idea that the connective tissue is physically coupled to
the muscle fibers which seems reasonable with regard to the available observations.
Nevertheless, additional experimental works must be done in order to support this assumption and to
study thoroughly the spatial organization of the myocardial collagen fibrils under normal and pathological
conditions.
\appendix
\setcounter{section}{0}
\section{Coordinate systems} \label{COORD}

\begin{figure}[H]
\label{coord}
\centering\includegraphics[scale=0.4]{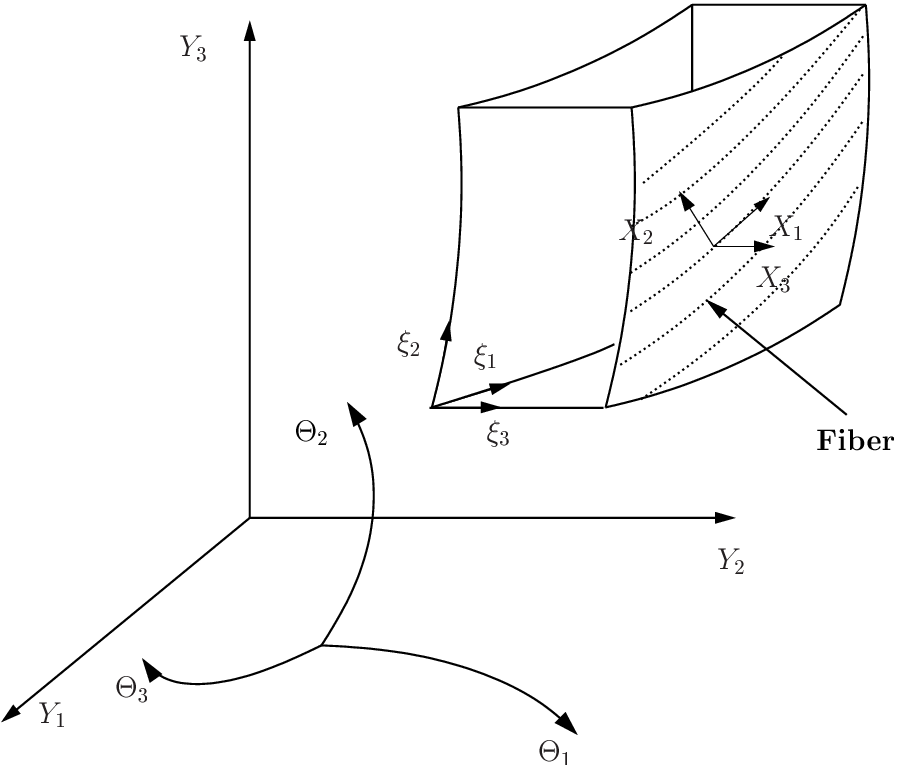}
\caption{Coordinate systems
(adapted from Costa et al. \cite{CHWWGM961}).}
\end{figure}

\begin{table}[H]
\caption{Notations for the coordinate systems used to formulate the finite element method
(adapted from Costa  et al. \cite{CHWWGM961}).
(I) Rectangular cartesian reference coordinates, (II) Curvilinear
world coordinates,   (III) Normalized finite element coordinate, (IV) Locally orthonormal body/fiber
coordinates (adapted from Costa  et al. \cite{CHWWGM961}).\label{tableau}}
\begin{center}
\begin{small}
$$\hspace{-1cm}\begin{array}{|c|c|c|c|c|c|c|c|}
\hline & \raisebox{-1.2ex}[0cm][0cm]{\hbox{\small
State}}&\raisebox{-1.2ex}[0cm][0cm]{ \hbox{\small
Indices}}&\raisebox{-1.2ex}[0cm][0cm]{\hbox{\small Coord.}} &
\hbox{\small Covariant}
&\hbox{\small Contravariant} &\multicolumn{2}{c|}{\raisebox{-1.2ex}[0cm][0cm]{\hbox{\small Metric tensors }}} \\
& & & &\hbox{\small  basis vectors} & \hbox{\small  basis vectors}&\multicolumn{2}{c|}{\hbox{}} \\
\hline
&P&R,S&Y^R&\e_R&\e_R&\delta_{RS}&\raisebox{0ex}[.5cm][0cm]{$\delta^{RS}$}\\
\cline{2-8}
\raisebox{2ex}[0.5cm][0.1cm]{(I)}&C&r,s&y^r&\e_r&\e_r&\delta_{rs}&\delta^{rs}\\
\hline
&P&A,B&\Theta^A&\raisebox{0ex}[1cm][.5cm]{$\ds\G_A^{(\theta)}=\frac{\partial{\mathbf{R}}}{\partial\Theta^A}$}&
\G^{(\theta)A}&G_{AB}^{(\theta)}&G^{(\theta)AB}\\
\cline{2-8} \raisebox{5ex}[0.5cm][0.1cm]{(II)}&C&\a,\b&\theta^\a&
\raisebox{0ex}[1cm][.5cm]{$\ds\mathbf{g}_\a^{(\theta)}=\frac{\partial\mathbf{r}}{\partial\theta^\a}$}&
\mathbf{g}^{(\theta)\a}&g_{\a\b}^{(\theta)}&g^{(\theta)\a\b}\\
\hline \rm
(III)&P&K,L&\xi^K&\raisebox{0ex}[1cm][.5cm]{$\ds\G_K^{(\xi)}=\frac{\partial\mathbf{R}}{\partial
\xi^K}$}&\raisebox{0ex}[.5cm][0cm]{$\G^{(\xi)K}$}&G_{KL}^{(\xi)}&G^{(\xi)KL}\\
\hline
&P&I,J&X^I&\raisebox{0ex}[1cm][.5cm]{$\ds\G_I^{(x)}=\frac{\partial\mathbf{R}}{\partial
X^I}$}&\raisebox{0ex}[.5cm][0cm]{$\G^{(x)I}$}&G_{IJ}^{(x)}=\d_{IJ}&G^{(x)IJ}=\d^{IJ}\\
&C&&&\raisebox{0ex}[1cm][.5cm]{$\ds\mathbf{g}_I^{(x)}=\frac{\partial\mathbf{r}}{\partial
X^I}$}&\ds\mathbf{g}^{(x)I}&g_{IJ}^{(x)}&g^{(x)IJ}\\
\cline{2-8}
\hline
\end{array}$$
\end{small}
\end{center}
\end{table}

\section{Second Piola-Kirchoff stress tensor at states $A$ and $C$}\label{P}

In the case of the free contraction state $A$, using
Eqs.(\ref{pa})-(\ref{CauchyC})-(\ref{taupa})-(\ref{Piola}),  we can write the components
$P_A^{IJ}$ of the second  Piola-Kirchoff stress tensor $\P_A$  in  the state $A$ base tensor
$(\mathbf{G}_I^{(x)} \otimes \mathbf{G}_J^{(x)})_A $ under the form:

\begin{eqnarray*}
P_A^{IJ}&=&-p_A \,g^{(x)IJ}+2G^{(x)IJ}\,W^*_1+(2\,W^*_4-q_A\,h_4) \,  f_P^{(x)I}
f_P^{(x)J}\nonumber\\
&&+\beta\,T^{(0)} f_A^{(x)I} f_A^{(x)J}-q_A\,h_6\, f_P^{\bot(x)I}
f_P^{\bot(x)J}
\end{eqnarray*}

 while in the case of the active loaded state $C$, using Eq.(\ref{CauchyA})-(\ref{Piola}), the components
$P_C^{IJ}$   of  $\P_C$ in the state $C$ base tensor $(\mathbf{G}_I^{(x)} \otimes
\mathbf{G}_J^{(x)})_C $  are

\begin{eqnarray*}
 P_C^{IJ}&=&-p_C   \,g^{(x)IJ}+2G^{(x)IJ}\,W^*_1+2\,W^*_4 \,  f_P^{(x)I}
f_P^{(x)J}\nonumber\\
&&+(\beta\,T^{(0)}+T_A^f)\, f_C^{(x)I} f_C^{(x)J}+T^{cf}_A\, f_C^{'(x)I}
f_C^{'(x)J}
\end{eqnarray*}

 \begin{eqnarray}\label{Wh}
{\rm where}\qquad  W^*_i&=&\frac{\partial W^*}{\partial I_i}= \frac{\partial W_{pas}}{\partial
I_i}+\beta \, \frac{\partial W_{act}^{f}}{\partial I_i}
+\delta_{AH}\, \frac{{\partial W_{\stackrel{pseudo}{ active}}}}{{\partial I_i }}
\quad i=1, \, 4\\
{\rm and}\qquad  h_i&=&\frac{\partial h(I_4,I_6)}{\partial I_i}\quad i=4, \, 6
 \end{eqnarray}

 $f_P^{(x)I}$, $f_P^{\bot(x)I} $   are respectively the components of the unit vectors $\f_P$,
$\f^\bot_P$    in the base ($\G_I^{(x)},\  I=1,2, 3$) and $f_H^{(x)I}$, $f_H^{'(x)I}$  with $H=A$ or
$C$, are respectively the components of the unit vectors $\f_H$, $\f'_H$ in the base
($\mathbf{g}_I^{(x)},\ I=1,2, 3$). The metric tensors $G^{(x)IJ}$, $g^{(x)IJ}$ are defined in table
1.\\
Following the definition of the locally orthonormal body/fiber coordinate system we have $
f_P^{(x)I}=\delta^{1I} $ and $ f_P^{\bot(x)I}=\delta^{2I}$.  On the other hand the
vectors $\f_H$ and $\f'_H$ are respectively defined through:

\begin{equation}
\f_H=\frac{\bPhi_{PH}\f_P}{\norm{\bPhi_{PH}\f_P}}=
\frac{\f^{(x)I}_P\mathbf{g}_I^{(x)}}{\norm{\f^{(x)I}_P\mathbf{g}_I^{(x)}}} \quad\hbox{and}\quad
\f'_H=\frac{\bPhi_{PH}\f^\bot_P}{\norm{\bPhi_{PH}\f^\bot_P}}=
\frac{\f^{\bot(x)I}_P\mathbf{g}_I^{(x)}}{\norm{\f^{\bot(x)I}_P\mathbf{g}_I^{(x)}}}
 \end{equation}

thus $\ds   f_H^{(x)I}=\frac{\delta^{1I}}{\norm{\mathbf{g}_1^{(x)}}}$ and $\ds
f_H^{'(x)I}=\frac{\delta^{2I}}{\norm{\mathbf{g}_2^{(x)}}}$ and we get finaly:

\begin{eqnarray}\label{PIJA}
 P_A^{IJ}&=&-p_A \,g^{(x)IJ}+2G^{(x)IJ}\,W^*_1+(2\,W^*_4-q_A\,h_4) \, \delta^{1I}\delta^{1J}
\nonumber\\ &&+\beta\,T^{(0)}
\frac{\delta^{1I}\delta^{1J}}{\norm{\mathbf{g}_1^{(x)}}^2}-q_A\,h_6\, \delta^{2I}\delta^{2J}
\end{eqnarray} and
\begin{eqnarray}\label{PIJC}
 P_C^{IJ}&=&-p_C   \,g^{(x)IJ}+2G^{(x)IJ}\,W^*_1+2\,W^*_4 \, \delta^{1I}\delta^{1J}\nonumber\\
&&+(\beta\,T^{(0)}+T_A^f)\, \frac{\delta^{1I}\delta^{1J}}{\norm{\mathbf{g}_1^{(x)}}^2}
+T^{cf}_A\, \frac{\delta^{2I}\delta^{2J}}{\norm{\mathbf{g}_2^{(x)}}^2}
\end{eqnarray}



\bibliographystyle{plain}

\end{document}